# HIGH-DIMENSIONAL GRAPHS AND VARIABLE SELECTION WITH THE LASSO

By Nicolai Meinshausen and Peter Bühlmann

*ETH Zürich*

The pattern of zero entries in the inverse covariance matrix of a multivariate normal distribution corresponds to conditional independence restrictions between variables. Covariance selection aims at estimating those structural zeros from data. We show that neighborhood selection with the Lasso is a computationally attractive alternative to standard covariance selection for sparse high-dimensional graphs. Neighborhood selection estimates the conditional independence restrictions separately for each node in the graph and is hence equivalent to variable selection for Gaussian linear models. We show that the proposed neighborhood selection scheme is consistent for sparse high-dimensional graphs. Consistency hinges on the choice of the penalty parameter. The oracle value for optimal prediction does not lead to a consistent neighborhood estimate. Controlling instead the probability of falsely joining some distinct connectivity components of the graph, consistent estimation for sparse graphs is achieved (with exponential rates), even when the number of variables grows as the number of observations raised to an arbitrary power.

**1. Introduction.** Consider the $p$-dimensional multivariate normal distributed random variable

$$X = (X_1, \ldots, X_p) \sim \mathcal{N}(\mu, \Sigma).$$

This includes Gaussian linear models where, for example, $X_1$ is the response variable and $\{X_k; 2 \le k \le p\}$ are the predictor variables. Assuming that the covariance matrix $\Sigma$ is nonsingular, the conditional independence structure of the distribution can be conveniently represented by a graphical model $\mathcal{G} = (\Gamma, E)$, where $\Gamma = \{1, \ldots, p\}$ is the set of nodes and $E$ the set of edges in $\Gamma \times \Gamma$. A pair $(a, b)$ is contained in the edge set $E$ if and only if $X_a$ is conditionally dependent on $X_b$, given all remaining variables $X_{\Gamma \setminus \{a,b\}} =$ $\{X_k; k \in \Gamma \setminus \{a, b\}\}$. Every pair of variables not contained in the edge set is conditionally independent, given all remaining variables, and corresponds to a zero entry in the inverse covariance matrix [12].

Covariance selection was introduced by Dempster [3] and aims at discovering the conditional independence restrictions (the graph) from a set of i.i.d. observations. Covariance selection traditionally relies on the discrete optimization of an objective function [5, 12]. Exhaustive search is computationally infeasible for all but very low-dimensional models. Usually, greedy forward or backward search is employed. In forward search, the initial estimate of the edge set is the empty set and edges are then added iteratively until a suitable stopping criterion is satisfied. The selection (deletion) of a single edge in this search strategy requires an MLE fit [15] for $O(p^2)$ different models. The procedure is not well suited for high-dimensional graphs. The existence of the MLE is not guaranteed in general if the number of observations is smaller than the number of nodes [1]. More disturbingly, the complexity of the procedure renders even greedy search strategies impractical for modestly sized graphs. In contrast, neighborhood selection with the Lasso, proposed in the following, relies on optimization of a convex function, applied consecutively to each node in the graph. The method is computationally very efficient and is consistent even for the high-dimensional setting, as will be shown.

Neighborhood selection is a subproblem of covariance selection. The neighborhood $\mathrm{ne}_a$ of a node $a \in \Gamma$ is the smallest subset of $\Gamma \setminus \{a\}$ so that, given all variables $X_{\mathrm{ne}_a}$ in the neighborhood, $X_a$ is conditionally independent of all remaining variables. The neighborhood of a node $a \in \Gamma$ consists of all nodes $b \in \Gamma \setminus \{a\}$ so that $(a, b) \in E$. Given $n$ i.i.d. observations of $X$, neighborhood selection aims at estimating (individually) the neighborhood of any given variable (or node). The neighborhood selection can be cast as a standard regression problem and can be solved efficiently with the Lasso [16], as will be shown in this paper.

The consistency of the proposed neighborhood selection will be shown for sparse high-dimensional graphs, where the number of variables is potentially growing as any power of the number of observations (high-dimensionality), whereas the number of neighbors of any variable is growing at most slightly slower than the number of observations (sparsity).

A number of studies have examined the case of regression with a growing number of parameters as sample size increases. The closest to our setting is the recent work of Greenshtein and Ritov [8], who study consistent prediction in a triangular setup very similar to ours (see also [10]). However, the problem of consistent estimation of the model structure, which is the relevant concept for graphical models, is very different and not treated in these studies.

We study in Section 2 under which conditions, and at which rate, the neighborhood estimate with the Lasso converges to the true neighborhood. The choice of the penalty is crucial in the high-dimensional setting. The oracle penalty for optimal prediction turns out to be inconsistent for estimation of the true model. This solution might include an unbounded number of noise variables in the model. We motivate a different choice of the penalty such that the probability of falsely connecting two or more distinct connectivity components of the graph is controlled at very low levels. Asymptotically, the probability of estimating the correct neighborhood converges exponentially to 1, even when the number of nodes in the graph is growing rapidly as any power of the number of observations. As a consequence, consistent estimation of the full edge set in a sparse high-dimensional graph is possible (Section 3).

Encouraging numerical results are provided in Section 4. The proposed estimate is shown to be both more accurate than the traditional forward selection MLE strategy and computationally much more efficient. The accuracy of the forward selection MLE fit is in particular poor if the number of nodes in the graph is comparable to the number of observations. In contrast, neighborhood selection with the Lasso is shown to be reasonably accurate for estimating graphs with several thousand nodes, using only a few hundred observations.

**2. Neighborhood selection.** Instead of assuming a fixed true underlying model, we adopt a more flexible approach similar to the triangular setup in [8]. Both the number of nodes in the graphs (number of variables), denoted by $p(n) = |\Gamma(n)|$, and the distribution (the covariance matrix) depend in general on the number of observations, so that $\Gamma = \Gamma(n)$ and $\Sigma = \Sigma(n)$. The neighborhood $\mathrm{ne}_a$ of a node $a \in \Gamma(n)$ is the smallest subset of $\Gamma(n) \setminus \{a\}$ so that $X_a$ is conditionally independent of all remaining variables. Denote the closure of node $a \in \Gamma(n)$ by $\mathrm{cl}_a := \mathrm{ne}_a \cup \{a\}$. Then

$$X_a \perp \{X_k; k \in \Gamma(n) \setminus \mathrm{cl}_a\} | X_{\mathrm{ne}_a}.$$

For details see [12]. The neighborhood depends in general on $n$ as well. However, this dependence is notationally suppressed in the following.

It is instructive to give a slightly different definition of a neighborhood. For each node $a \in \Gamma(n)$, consider optimal prediction of $X_a$, given all remaining variables. Let $\theta^a \in \mathbb{R}^{p(n)}$ be the vector of coefficients for optimal prediction,

$$(1) \qquad \theta^a = \underset{\theta : \theta_a = 0}{\arg\min} E\left(X_a - \sum_{k \in \Gamma(n)} \theta_k X_k\right)^2.$$

As a generalization of (1), which will be of use later, consider optimal prediction of $X_a$, given only a subset of variables $\{X_k; k \in \mathcal{A}\}$, where $\mathcal{A} \subseteq$

$\Gamma(n) \setminus \{a\}$. The optimal prediction is characterized by the vector $\theta^{a,\mathcal{A}}$,

$$
(2) \qquad \theta^{a,\mathcal{A}} = \underset{\theta \,:\, \theta_k = 0, \forall k \notin \mathcal{A}}{\arg\min} E\left(X_a - \sum_{k \in \Gamma(n)} \theta_k X_k\right)^2.
$$

The elements of $\theta^a$ are determined by the inverse covariance matrix [12]. For $b \in \Gamma \setminus \{a\}$ and $K(n) = \Sigma^{-1}(n)$, it holds that $\theta_b^a = -K_{ab}(n)/K_{aa}(n)$. The set of nonzero coefficients of $\theta^a$ is identical to the set $\{b \in \Gamma(n) \setminus \{a\} : K_{ab}(n) \neq 0\}$ of nonzero entries in the corresponding row vector of the inverse covariance matrix and defines precisely the set of neighbors of node $a$. The best predictor for $X_a$ is thus a linear function of variables in the set of neighbors of the node $a$ only. The set of neighbors of a node $a \in \Gamma(n)$ can hence be written as

$$\mathrm{ne}_a = \{b \in \Gamma(n) : \theta_b^a \neq 0\}.$$

This set corresponds to the set of effective predictor variables in regression with response variable $X_a$ and predictor variables $\{X_k; k \in \Gamma(n) \setminus \{a\}\}$. Given $n$ independent observations of $X \sim \mathcal{N}(0, \Sigma(n))$, neighborhood selection tries to estimate the set of neighbors of a node $a \in \Gamma(n)$. As the optimal linear prediction of $X_a$ has nonzero coefficients precisely for variables in the set of neighbors of the node $a$, it seems reasonable to try to exploit this relation.

2.1. *Neighborhood selection with the Lasso.* It is well known that the Lasso, introduced by Tibshirani [16], and known as Basis Pursuit in the context of wavelet regression [2], has a parsimonious property [11]. When predicting a variable $X_a$ with all remaining variables $\{X_k; k \in \Gamma(n) \setminus \{a\}\}$, the vanishing Lasso coefficient estimates identify asymptotically the neighborhood of node $a$ in the graph, as shown in the following. Let the $n \times p(n)$-dimensional matrix $\mathbf{X}$ contain $n$ independent observations of $X$, so that the columns $\mathbf{X}_a$ correspond for all $a \in \Gamma(n)$ to the vector of $n$ independent observations of $X_a$. Let $\langle \cdot, \cdot \rangle$ be the usual inner product on $\mathbb{R}^n$ and $\| \cdot \|_2$ the corresponding norm.

The Lasso estimate $\hat{\theta}^{a,\lambda}$ of $\theta^a$ is given by

$$
(3) \qquad \hat{\theta}^{a,\lambda} = \underset{\theta \,:\, \theta_a = 0}{\arg\min}(n^{-1}\|\mathbf{X}_a - \mathbf{X}\theta\|_2^2 + \lambda\|\theta\|_1),
$$

where $\|\theta\|_1 = \sum_{b \in \Gamma(n)} |\theta_b|$ is the $l_1$-norm of the coefficient vector. Normalization of all variables to a common empirical variance is recommended for the estimator in (3). The solution to (3) is not necessarily unique. However, if uniqueness fails, the set of solutions is still convex and all our results about neighborhoods (in particular Theorems 1 and 2) hold for any solution of (3).

Other regression estimates have been proposed, which are based on the $l_p$-norm, where $p$ is typically in the range $[0, 2]$ (see [7]). A value of $p = 2$ leads to the ridge estimate, while $p = 0$ corresponds to traditional model selection. It is well known that the estimates have a parsimonious property (with some components being exactly zero) for $p \leq 1$ only, while the optimization problem in (3) is only convex for $p \geq 1$. Hence $l_1$-constrained empirical risk minimization occupies a unique position, as $p = 1$ is the only value of $p$ for which variable selection takes place while the optimization problem is still convex and hence feasible for high-dimensional problems.

The neighborhood estimate (parameterized by $\lambda$) is defined by the nonzero coefficient estimates of the $l_1$-penalized regression,

$$\hat{\text{ne}}_a^\lambda = \{b \in \Gamma(n) : \hat{\theta}_b^{a,\lambda} \neq 0\}.$$

Each choice of a penalty parameter $\lambda$ specifies thus an estimate of the neighborhood $\text{ne}_a$ of node $a \in \Gamma(n)$ and one is left with the choice of a suitable penalty parameter. Larger values of the penalty tend to shrink the size of the estimated set, while more variables are in general included into $\hat{\text{ne}}_a^\lambda$ if the value of $\lambda$ is diminished.

2.2. *The prediction-oracle solution.* A seemingly useful choice of the penalty parameter is the (unavailable) prediction-oracle value,

$$\lambda_{\text{oracle}} = \arg\min_\lambda E\left(X_a - \sum_{k \in \Gamma(n)} \hat{\theta}_k^{a,\lambda} X_k\right)^2.$$

The expectation is understood to be with respect to a new $X$, which is independent of the sample on which $\hat{\theta}^{a,\lambda}$ is estimated. The prediction-oracle penalty minimizes the predictive risk among all Lasso estimates. An estimate of $\lambda_{\text{oracle}}$ is obtained by the cross-validated choice $\lambda_{\text{cv}}$.

For $l_0$-penalized regression it was shown by Shao [14] that the cross-validated choice of the penalty parameter is consistent for model selection under certain conditions on the size of the validation set. The prediction-oracle solution does not lead to consistent model selection for the Lasso, as shown in the following for a simple example.

PROPOSITION 1. *Let the number of variables grow to infinity, $p(n) \to \infty$, for $n \to \infty$, with $p(n) = o(n^\gamma)$ for some $\gamma > 0$. Assume that the covariance matrices $\Sigma(n)$ are identical to the identity matrix except for some pair $(a, b) \in \Gamma(n) \times \Gamma(n)$, for which $\Sigma_{ab}(n) = \Sigma_{ba}(n) = s$, for some $0 < s < 1$ and all $n \in \mathbb{N}$. The probability of selecting the wrong neighborhood for node $a$ converges to 1 under the prediction-oracle penalty,*

$$P(\hat{\text{ne}}_a^{\lambda_{\text{oracle}}} \neq \text{ne}_a) \to 1 \quad \text{for } n \to \infty.$$

A proof is given in the Appendix. It follows from the proof of Proposition 1 that many noise variables are included in the neighborhood estimate with the prediction-oracle solution. In fact, the probability of including noise variables with the prediction-oracle solution does not even vanish asymptotically for a fixed number of variables. If the penalty is chosen larger than the prediction-optimal value, consistent neighborhood selection is possible with the Lasso, as demonstrated in the following.

2.3. *Assumptions.* We make a few assumptions to prove consistency of neighborhood selection with the Lasso. We always assume availability of $n$ independent observations from $X \sim \mathcal{N}(0, \Sigma)$.

*High-dimensionality.* The number of variables is allowed to grow as the number of observations $n$ raised to an arbitrarily high power.

ASSUMPTION 1. There exists $\gamma > 0$, so that

$$p(n) = O(n^\gamma) \qquad \text{for } n \to \infty.$$

In particular, it is allowed for the following analysis that the number of variables is very much larger than the number of observations, $p(n) \gg n$.

*Nonsingularity.* We make two regularity assumptions for the covariance matrices.

ASSUMPTION 2. For all $a \in \Gamma(n)$ and $n \in \mathbb{N}$, $\text{Var}(X_a) = 1$. There exists $v^2 > 0$, so that for all $n \in \mathbb{N}$ and $a \in \Gamma(n)$,

$$\text{Var}(X_a | X_{\Gamma(n) \setminus \{a\}}) \geq v^2.$$

Common variance can always be achieved by appropriate scaling of the variables. A scaling to a common (empirical) variance of all variables is desirable, as the solutions would otherwise depend on the chosen units or dimensions in which they are represented. The second part of the assumption explicitly excludes singular or nearly singular covariance matrices. For singular covariance matrices, edges are not uniquely defined by the distribution and it is hence not surprising that nearly singular covariance matrices are not suitable for consistent variable selection. Note, however, that the *empirical* covariance matrix is a.s. singular if $p(n) > n$, which is allowed in our analysis.

*Sparsity.* The main assumption is the sparsity of the graph. This entails a restriction on the size of the neighborhoods of variables.

ASSUMPTION 3. There exists some $0 \le \kappa < 1$ so that

$$\max_{a \in \Gamma(n)} |\mathrm{ne}_a| = O(n^\kappa) \qquad \text{for } n \to \infty.$$

This assumption limits the maximal possible rate of growth for the size of neighborhoods.

For the next sparsity condition, consider again the definition in (2) of the optimal coefficient $\theta^{b,\mathcal{A}}$ for prediction of $X_b$, given variables in the set $\mathcal{A} \subset \Gamma(n)$.

ASSUMPTION 4. There exists some $\vartheta < \infty$ so that for all neighboring nodes $a, b \in \Gamma(n)$ and all $n \in \mathbb{N}$,

$$\|\theta^{a,\mathrm{ne}_b \setminus \{a\}}\|_1 \le \vartheta.$$

This assumption is, for example, satisfied if Assumption 2 holds and the size of the overlap of neighborhoods is bounded by an arbitrarily large number from above for neighboring nodes. That is, if there exists some $m < \infty$ so that for all $n \in \mathbb{N}$,

(4) $$\max_{a,b \in \Gamma(n), b \in \mathrm{ne}_a} |\mathrm{ne}_a \cap \mathrm{ne}_b| \le m \qquad \text{for } n \to \infty,$$

then Assumption 4 is satisfied. To see this, note that Assumption 2 gives a finite bound for the $l_2$-norm of $\theta^{a,\mathrm{ne}_b \setminus \{a\}}$, while (4) gives a finite bound for the $l_0$-norm. Taken together, Assumption 4 is implied.

*Magnitude of partial correlations.* The next assumption bounds the magnitude of partial correlations from below. The partial correlation $\pi_{ab}$ between variables $X_a$ and $X_b$ is the correlation after having eliminated the linear effects from all remaining variables $\{X_k; k \in \Gamma(n) \setminus \{a,b\}\}$; for details see [12].

ASSUMPTION 5. There exist a constant $\delta > 0$ and some $\xi > \kappa$, with $\kappa$ as in Assumption 3, so that for every $(a,b) \in E$,

$$|\pi_{ab}| \ge \delta n^{-(1-\xi)/2}.$$

It will be shown below that Assumption 5 cannot be relaxed in general. Note that neighborhood selection for node $a \in \Gamma(n)$ is equivalent to simultaneously testing the null hypothesis of zero partial correlation between variable $X_a$ and all remaining variables $X_b$, $b \in \Gamma(n) \setminus \{a\}$. The null hypothesis of zero partial correlation between two variables can be tested by using the corresponding entry in the normalized inverse empirical covariance matrix. A graph estimate based on such tests has been proposed by Drton and Perlman [4]. Such a test can only be applied, however, if the number of variables is smaller than the number of observations, $p(n) \leq n$, as the empirical covariance matrix is singular otherwise. Even if $p(n) = n - c$ for some constant $c > 0$, Assumption 5 would have to hold with $\xi = 1$ to have a positive power of rejecting false null hypotheses for such an estimate; that is, partial correlations would have to be bounded by a positive value from below.

*Neighborhood stability.* The last assumption is referred to as neighborhood stability. Using the definition of $\theta^{a,\mathcal{A}}$ in (2), define for all $a, b \in \Gamma(n)$,

$$(5) \qquad S_a(b) := \sum_{k \in \mathrm{ne}_a} \mathrm{sign}(\theta_k^{a,\mathrm{ne}_a}) \theta_k^{b,\mathrm{ne}_a}.$$

The assumption of neighborhood stability restricts the magnitude of the quantities $S_a(b)$ for nonneighboring nodes $a, b \in \Gamma(n)$.

ASSUMPTION 6. There exists some $\delta < 1$ so that for all $a, b \in \Gamma(n)$ with $b \notin \mathrm{ne}_a$,

$$|S_a(b)| < \delta.$$

It is shown in Proposition 3 that this assumption cannot be relaxed.

We give in the following a more intuitive condition which essentially implies Assumption 6. This will justify the term neighborhood stability. Consider the definition in (1) of the optimal coefficients $\theta^a$ for prediction of $X_a$. For $\eta > 0$, define $\theta^a(\eta)$ as the optimal set of coefficients under an additional $l_1$-penalty,

$$(6) \qquad \theta^a(\eta) := \underset{\theta \,:\, \theta_a = 0}{\arg\min} E\left(X_a - \sum_{k \in \Gamma(n)} \theta_k X_k\right)^2 + \eta \|\theta\|_1.$$

The neighborhood $\mathrm{ne}_a$ of node $a$ was defined as the set of nonzero coefficients of $\theta^a$, $\mathrm{ne}_a = \{k \in \Gamma(n) : \theta_k^a \neq 0\}$. Define the disturbed neighborhood $\mathrm{ne}_a(\eta)$ as

$$\mathrm{ne}_a(\eta) := \{k \in \Gamma(n) : \theta_k^a(\eta) \neq 0\}.$$

It clearly holds that $\mathrm{ne}_a = \mathrm{ne}_a(0)$. The assumption of neighborhood stability is satisfied if there exists some infinitesimally small perturbation $\eta$, which may depend on $n$, so that the disturbed neighborhood $\mathrm{ne}_a(\eta)$ is identical to the undisturbed neighborhood $\mathrm{ne}_a(0)$.

PROPOSITION 2. *If there exists some $\eta > 0$ so that $\mathrm{ne}_a(\eta) = \mathrm{ne}_a(0)$, then $|S_a(b)| \leq 1$ for all $b \in \Gamma(n) \setminus \mathrm{ne}_a$.*

A proof is given in the Appendix.

In light of Proposition 2 it seems that Assumption 6 is a very weak condition. To give one example, Assumption 6 is automatically satisfied under the much stronger assumption that the graph does not contain cycles. We give a brief reasoning for this. Consider two nonneighboring nodes $a$ and $b$. If the nodes are in different connectivity components, there is nothing left to show as $S_a(b) = 0$. If they are in the same connectivity component, then there exists one node $k \in \mathrm{ne}_a$ that separates $b$ from $\mathrm{ne}_a \setminus \{k\}$, as there is just one unique path between any two variables in the same connectivity component if the graph does not contain cycles. Using the global Markov property, the random variable $X_b$ is independent of $X_{\mathrm{ne}_a \setminus \{k\}}$, given $X_k$. The random variable $E(X_b|X_{\mathrm{ne}_a})$ is thus a function of $X_k$ only. As the distribution is Gaussian, $E(X_b|X_{\mathrm{ne}_a}) = \theta_k^{b,\mathrm{ne}_a} X_k$. By Assumption 2, $\mathrm{Var}(X_b|X_{\mathrm{ne}_a}) = v^2$ for some $v^2 > 0$. It follows that $\mathrm{Var}(X_b) = v^2 + (\theta_k^{b,\mathrm{ne}_a})^2 = 1$ and hence $\theta_k^{b,\mathrm{ne}_a} = \sqrt{1-v^2} < 1$, which implies that Assumption 6 is indeed satisfied if the graph does not contain cycles.

We mention that Assumption 6 is likewise satisfied if the inverse covariance matrices $\Sigma^{-1}(n)$ are for each $n \in \mathbb{N}$ diagonally dominant. A matrix is said to be diagonally dominant if and only if, for each row, the sum of the absolute values of the nondiagonal elements is smaller than the absolute value of the diagonal element. The proof of this is straightforward but tedious and hence is omitted.

2.4. *Controlling type I errors.* The asymptotic properties of Lasso-type estimates in regression have been studied in detail by Knight and Fu [11] for a fixed number of variables. Their results say that the penalty parameter $\lambda$ should decay for an increasing number of observations at least as fast as $n^{-1/2}$ to obtain an $n^{1/2}$-consistent estimate. It turns out that a slower rate is needed for consistent model selection in the high-dimensional case where $p(n) \gg n$. However, a rate $n^{-(1-\varepsilon)/2}$ with any $\kappa < \varepsilon < \xi$ (where $\kappa, \xi$ are defined as in Assumptions 3 and 5) is sufficient for consistent neighborhood selection, even when the number of variables is growing rapidly with the number of observations.

THEOREM 1. *Let Assumptions 1–6 hold. Let the penalty parameter satisfy $\lambda_n \sim dn^{-(1-\varepsilon)/2}$ with some $\kappa < \varepsilon < \xi$ and $d > 0$. There exists some $c > 0$ so that, for all $a \in \Gamma(n)$,*

$$P(\hat{\mathrm{ne}}_a^\lambda \subseteq \mathrm{ne}_a) = 1 - O(\exp(-cn^\varepsilon)) \qquad \text{for } n \to \infty.$$

A proof is given in the Appendix.

Theorem 1 states that the probability of (falsely) including any of the nonneighboring variables of the node $a \in \Gamma(n)$ into the neighborhood estimate vanishes exponentially fast, even though the number of nonneighboring variables may grow very rapidly with the number of observations. It is shown in the following that Assumption 6 cannot be relaxed.

PROPOSITION 3. *If there exists some $a, b \in \Gamma(n)$ with $b \notin \mathrm{ne}_a$ and $|S_a(b)| > 1$, then, for $\lambda = \lambda_n$ as in Theorem 1,*

$$P(\hat{\mathrm{ne}}_a^\lambda \subseteq \mathrm{ne}_a) \to 0 \qquad \text{for } n \to \infty.$$

A proof is given in the Appendix. Assumption 6 of neighborhood stability is hence critical for the success of Lasso neighborhood selection.

2.5. *Controlling type* II *errors.* So far it has been shown that the probability of falsely including variables into the neighborhood can be controlled by the Lasso. The question arises whether the probability of including all neighboring variables into the neighborhood estimate converges to 1 for $n \to \infty$.

THEOREM 2. *Let the assumptions of Theorem 1 be satisfied. For $\lambda = \lambda_n$ as in Theorem 1, for some $c > 0$*

$$P(\mathrm{ne}_a \subseteq \hat{\mathrm{ne}}_a^\lambda) = 1 - O(\exp(-cn^\varepsilon)) \qquad \text{for } n \to \infty.$$

A proof is given in the Appendix.

It may be of interest whether Assumption 5 could be relaxed, so that edges are detected even if the partial correlation is vanishing at a rate $n^{-(1-\xi)/2}$ for some $\xi < \kappa$. The following proposition says that $\xi > \varepsilon$ (and thus $\xi > \kappa$ as $\varepsilon > \kappa$) is a necessary condition if a stronger version of Assumption 4 holds, which is satisfied for forests and trees, for example.

PROPOSITION 4. *Let the assumptions of Theorem 1 hold with $\vartheta < 1$ in Assumption 4, except that for $a \in \Gamma(n)$, let there be some $b \in \Gamma(n) \setminus \{a\}$ with $\pi_{ab} \neq 0$ and $|\pi_{ab}| = O(n^{-(1-\xi)/2})$ for $n \to \infty$ for some $\xi < \varepsilon$. Then*

$$P(b \in \hat{\mathrm{ne}}_a^\lambda) \to 0 \qquad \text{for } n \to \infty.$$

Theorem 2 and Proposition 4 say that edges between nodes for which partial correlation vanishes at a rate $n^{-(1-\xi)/2}$ are, with probability converging to 1 for $n \to \infty$, detected if $\xi > \varepsilon$ and are undetected if $\xi < \varepsilon$. The results do not cover the case $\xi = \varepsilon$, which remains a challenging question for further research.

All results so far have treated the distinction between zero and nonzero partial correlations only. The signs of partial correlations of neighboring nodes can be estimated consistently under the same assumptions and with the same rates, as can be seen in the proofs.

**3. Covariance selection.** It follows from Section 2 that it is possible under certain conditions to estimate the neighborhood of each node in the graph consistently, for example,

$$P(\hat{\text{ne}}_a^\lambda = \text{ne}_a) \to 1 \quad \text{for } n \to \infty.$$

The full graph is given by the set $\Gamma(n)$ of nodes and the edge set $E = E(n)$. The edge set contains those pairs $(a,b) \in \Gamma(n) \times \Gamma(n)$ for which the partial correlation between $X_a$ and $X_b$ is not zero. As the partial correlations are precisely nonzero for neighbors, the edge set $E \subseteq \Gamma(n) \times \Gamma(n)$ is given by

$$E = \{(a,b) : a \in \text{ne}_b \land b \in \text{ne}_a\}.$$

The first condition, $a \in \text{ne}_b$, implies in fact the second, $b \in \text{ne}_a$, and vice versa, so that the edge is as well identical to $\{(a,b) : a \in \text{ne}_b \lor b \in \text{ne}_a\}$. For an estimate of the edge set of a graph, we can apply neighborhood selection to each node in the graph. A natural estimate of the edge set is then given by $\hat{E}^{\lambda,\land} \subseteq \Gamma(n) \times \Gamma(n)$, where

(7) $$\hat{E}^{\lambda,\land} = \{(a,b) : a \in \hat{\text{ne}}_b^\lambda \land b \in \hat{\text{ne}}_a^\lambda\}.$$

Note that $a \in \hat{\text{ne}}_b^\lambda$ does not necessarily imply $b \in \hat{\text{ne}}_a^\lambda$ and vice versa. We can hence also define a second, less conservative, estimate of the edge set by

(8) $$\hat{E}^{\lambda,\lor} = \{(a,b) : a \in \hat{\text{ne}}_b^\lambda \lor b \in \hat{\text{ne}}_a^\lambda\}.$$

The discrepancies between the estimates (7) and (8) are quite small in our experience. Asymptotically the difference between both estimates vanishes, as seen in the following corollary. We refer to both edge set estimates collectively with the generic notation $\hat{E}^\lambda$, as the following result holds for both of them.

COROLLARY 1. *Under the conditions of Theorem 2, for some $c > 0$,*

$$P(\hat{E}^\lambda = E) = 1 - O(\exp(-cn^\varepsilon)) \quad \text{for } n \to \infty.$$

The claim follows since $|\Gamma(n)|^2 = p(n)^2 = O(n^{2\gamma})$ by Assumption 1 and neighborhood selection has an exponentially fast convergence rate as described by Theorem 2. Corollary 1 says that the conditional independence structure of a multivariate normal distribution can be estimated consistently by combining the neighborhood estimates for all variables.

Note that there are in total $2^{(p^2-p)/2}$ distinct graphs for a $p$-dimensional variable. However, for each of the $p$ nodes there are only $2^{p-1}$ distinct potential neighborhoods. By breaking the graph selection problem into a consecutive series of neighborhood selection problems, the complexity of the search is thus reduced substantially at the price of potential inconsistencies between neighborhood estimates. Graph estimates that apply this strategy for complexity reduction are sometimes called dependency networks [9]. The complexity of the proposed neighborhood selection for one node with the Lasso is reduced further to $O(np \min\{n,p\})$, as the Lars procedure of Efron, Hastie, Johnstone and Tibshirani [6] requires $O(\min\{n,p\})$ steps, each of complexity $O(np)$. For high-dimensional problems as in Theorems 1 and 2, where the number of variables grows as $p(n) \sim cn^\gamma$ for some $c > 0$ and $\gamma > 1$, this is equivalent to $O(p^{2+2/\gamma})$ computations for the whole graph. The complexity of the proposed method thus scales approximately quadratic with the number of nodes for large values of $\gamma$.

Before providing some numerical results, we discuss in the following the choice of the penalty parameter.

*Finite-sample results and significance.* It was shown above that consistent neighborhood and covariance selection is possible with the Lasso in a high-dimensional setting. However, the asymptotic considerations give little advice on how to choose a specific penalty parameter for a given problem. Ideally, one would like to guarantee that pairs of variables which are not contained in the edge set enter the estimate of the edge set only with very low (prespecified) probability. Unfortunately, it seems very difficult to obtain such a result as the probability of falsely including a pair of variables into the estimate of the edge set depends on the exact covariance matrix, which is in general unknown. It is possible, however, to constrain the probability of (falsely) connecting two distinct connectivity components of the true graph. The connectivity component $C_a \subseteq \Gamma(n)$ of a node $a \in \Gamma(n)$ is the set of nodes which are connected to node $a$ by a chain of edges. The neighborhood $\text{ne}_a$ is clearly part of the connectivity component $C_a$.

Let $\hat{C}_a^\lambda$ be the connectivity component of $a$ in the estimated graph $(\Gamma, \hat{E}^\lambda)$. For any level $0 < \alpha < 1$, consider the choice of the penalty

$$\lambda(\alpha) = \frac{2\hat{\sigma}_a}{\sqrt{n}} \tilde{\Phi}^{-1}\left(\frac{\alpha}{2p(n)^2}\right), \tag{9}$$

where $\tilde{\Phi} = 1 - \Phi$ [$\Phi$ is the c.d.f. of $\mathcal{N}(0,1)$] and $\hat{\sigma}_a^2 = n^{-1}\langle \mathbf{X}_a, \mathbf{X}_a \rangle$. The probability of falsely joining two distinct connectivity components with the estimate of the edge set is bounded by the level $\alpha$ under the choice $\lambda = \lambda(\alpha)$ of the penalty parameter, as shown in the following theorem.

THEOREM 3. *Let Assumptions* 1–6 *be satisfied. Using the penalty parameter* $\lambda(\alpha)$, *we have for all* $n \in \mathbb{N}$ *that*

$$P(\exists\, a \in \Gamma(n) : \hat{C}_a^\lambda \nsubseteq C_a) \leq \alpha.$$

A proof is given in the Appendix. This implies that if the edge set is empty ($E = \varnothing$), it is estimated by an empty set with high probability,

$$P(\hat{E}^\lambda = \varnothing) \geq 1 - \alpha.$$

Theorem 3 is a finite-sample result. The previous asymptotic results in Theorems 1 and 2 hold if the level $\alpha$ vanishes exponentially to zero for an increasing number of observations, leading to consistent edge set estimation.

**4. Numerical examples.** We use both the Lasso estimate from Section 3 and forward selection MLE [5, 12] to estimate sparse graphs. We found it difficult to compare numerically neighborhood selection with forward selection MLE for more than 30 nodes in the graph. The high computational complexity of the forward selection MLE made the computations for such relatively low-dimensional problems very costly already. The Lasso scheme in contrast handled with ease graphs with more than 1000 nodes, using the recent algorithm developed in [6]. Where comparison was feasible, the performance of the neighborhood selection scheme was better. The difference was particularly pronounced if the ratio of observations to variables was low, as can be seen in Table 1, which will be described in more detail below.

First we give an account of the generation of the underlying graphs which we are trying to estimate. A realization of an underlying (random) graph is given in the left panel of Figure 1. The nodes of the graph are associated with spatial location and the location of each node is distributed identically and uniformly in the two-dimensional square $[0,1]^2$. Every pair of nodes is included initially in the edge set with probability $\varphi(d/\sqrt{p})$, where $d$ is the

TABLE 1
*The average number of correctly identified edges as a function of the number $k$ of falsely included edges for $n = 40$ observations and $p = 10, 20, 30$ nodes for forward selection MLE (FS), $\hat{E}^{\lambda,\vee}$, $\hat{E}^{\lambda,\wedge}$ and random guessing*

|  | $p = 10$ | | | $p = 20$ | | | $p = 30$ | | |
|---|---|---|---|---|---|---|---|---|---|
| $k$ | 0 | 5 | 10 | 0 | 5 | 10 | 0 | 5 | 10 |
| Random | 0.2 | 1.9 | 3.7 | 0.1 | 0.7 | 1.4 | 0.1 | 0.5 | 0.9 |
| FS | 7.6 | 14.1 | 17.1 | 8.9 | 16.6 | 21.6 | 0.6 | 1.8 | 3.2 |
| $\hat{E}^{\lambda,\vee}$ | 8.2 | 15.0 | 17.6 | 9.3 | 18.5 | 23.9 | 11.4 | 21.4 | 26.3 |
| $\hat{E}^{\lambda,\wedge}$ | 8.5 | 14.7 | 17.6 | 9.5 | 19.1 | 34.0 | 14.1 | 21.4 | 27.4 |

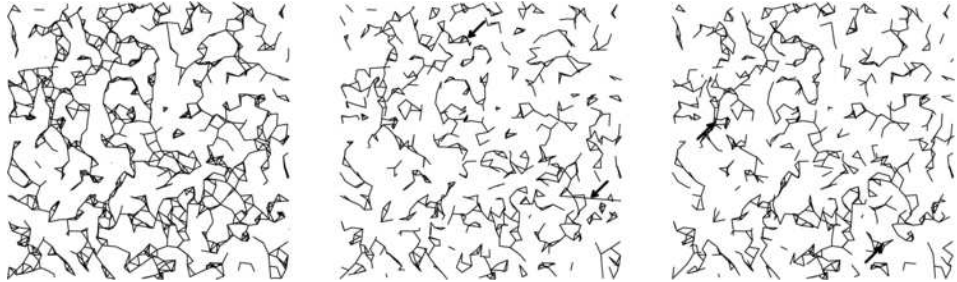

FIG. 1. *A realization of a graph is shown on the left, generated as described in the text. The graph consists of* 1000 *nodes and* 1747 *edges out of* 449,500 *distinct pairs of variables. The estimated edge set, using estimate* (7) *at level* $\alpha = 0.05$ [*see* (9)], *is shown in the middle. There are two erroneously included edges, marked by an arrow, while* 1109 *edges are correctly detected. For estimate* (8) *and an adjusted level as described in the text, the result is shown on the right. Again two edges are erroneously included. Not a single pair of disjoint connectivity components of the true graph has been (falsely) joined by either estimate.*

Euclidean distance between the pair of variables and $\varphi$ is the density of the standard normal distribution. The maximum number of edges connecting to each node is limited to four to achieve the desired sparsity of the graph. Edges which connect to nodes which do not satisfy this constraint are removed randomly until the constraint is satisfied for all edges. Initially all variables have identical conditional variance and the partial correlation between neighbors is set to 0.245 (absolute values less than 0.25 guarantee positive definiteness of the inverse covariance matrix); that is, $\Sigma_{aa}^{-1} = 1$ for all nodes $a \in \Gamma$, $\Sigma_{ab}^{-1} = 0.245$ if there is an edge connecting $a$ and $b$ and $\Sigma_{ab}^{-1} = 0$ otherwise. The diagonal elements of the corresponding covariance matrix are in general larger than 1. To achieve constant variance, all variables are finally rescaled so that the diagonal elements of $\Sigma$ are all unity. Using the Cholesky transformation of the covariance matrix, $n$ independent samples are drawn from the corresponding Gaussian distribution.

The average number of edges which are correctly included into the estimate of the edge set is shown in Table 1 as a function of the number of edges which are falsely included. The accuracy of the forward selection MLE is comparable to the proposed Lasso neighborhood selection if the number of nodes is much smaller than the number of observations. The accuracy of the forward selection MLE breaks down, however, if the number of nodes is approximately equal to the number of observations. Forward selection MLE is only marginally better than random guessing in this case. Computation of the forward selection MLE (using MIM, [5]) on the same desktop took up to several hundred times longer than the Lasso neighborhood selection for the full graph. For more than 30 nodes, the differences are even more pronounced.

The Lasso neighborhood selection can be applied to hundred- or thousand-dimensional graphs, a realistic size, for example, biological networks. A graph with 1000 nodes (following the same model as described above) and its estimates (7) and (8), using 600 observations, are shown in Figure 1. A level $\alpha = 0.05$ is used for the estimate $\hat{E}^{\lambda,\vee}$. For better comparison, the level $\alpha$ was adjusted to $\alpha = 0.064$ for the estimate $\hat{E}^{\lambda,\wedge}$, so that both estimates lead to the same number of included edges. There are two erroneous edge inclusions, while 1109 out of all 1747 edges have been correctly identified by either estimate. Of these 1109 edges, 907 are common to both estimates while 202 are just present in either (7) or (8).

To examine if results are critically dependent on the assumption of Gaussianity, long-tailed noise is added to the observations. Instead of $n$ i.i.d. observations of $X \sim \mathcal{N}(0, \Sigma)$, $n$ i.i.d. observations of $X + 0.1Z$ are made, where the components of $Z$ are independent and follow a $t_2$-distribution. For 10 simulations (with each 500 observations), the proportion of false rejections among all rejections increases only slightly from 0.8% (without long-tailed noise) to 1.4% (with long tailed-noise) for $\hat{E}^{\lambda,\vee}$ and from 4.8% to 5.2% for $\hat{E}^{\lambda,\wedge}$. Our limited numerical experience suggests that the properties of the graph estimator do not seem to be critically affected by deviations from Gaussianity.

## APPENDIX: PROOFS

**A.1. Notation and useful lemmas.** As a generalization of (3), the Lasso estimate $\hat{\theta}^{a,\mathcal{A},\lambda}$ of $\theta^{a,\mathcal{A}}$, defined in (2), is given by

$$(A.1) \qquad \hat{\theta}^{a,\mathcal{A},\lambda} = \operatorname*{arg\,min}_{\theta\,:\,\theta_k = 0\,\forall\,k \notin \mathcal{A}} (n^{-1}\|\mathbf{X}_a - \mathbf{X}\theta\|_2^2 + \lambda\|\theta\|_1).$$

The notation $\hat{\theta}^{a,\lambda}$ is thus just a shorthand notation for $\hat{\theta}^{a,\Gamma(n)\setminus\{a\},\lambda}$.

LEMMA A.1. *Given $\theta \in \mathbb{R}^{p(n)}$, let $G(\theta)$ be a $p(n)$-dimensional vector with elements*

$$G_b(\theta) = -2n^{-1}\langle \mathbf{X}_a - \mathbf{X}\theta, \mathbf{X}_b \rangle.$$

*A vector $\hat{\theta}$ with $\hat{\theta}_k = 0, \forall\,k \in \Gamma(n) \setminus \mathcal{A}$ is a solution to (A.1) iff for all $b \in \mathcal{A}$, $G_b(\hat{\theta}) = -\operatorname{sign}(\hat{\theta}_b)\lambda$ in case $\hat{\theta}_b \neq 0$ and $|G_b(\hat{\theta})| \leq \lambda$ in case $\hat{\theta}_b = 0$. Moreover, if the solution is not unique and $|G_b(\hat{\theta})| < \lambda$ for some solution $\hat{\theta}$, then $\hat{\theta}_b = 0$ for all solutions of (A.1).*

PROOF. Denote the subdifferential of

$$n^{-1}\|\mathbf{X}_a - \mathbf{X}\theta\|_2^2 + \lambda\|\theta\|_1$$

with respect to $\theta$ by $D(\theta)$. The vector $\hat{\theta}$ is a solution to (A.1) iff there exists an element $d \in D(\hat{\theta})$ so that $d_b = 0$, $\forall b \in \mathcal{A}$. $D(\theta)$ is given by $\{G(\theta) + \lambda e, e \in S\}$, where $S \subset \mathbb{R}^{p(n)}$ is given by $S := \{e \in \mathbb{R}^{p(n)} : e_b = \text{sign}(\theta_b) \text{ if } \theta_b \neq 0 \text{ and } e_b \in [-1,1] \text{ if } \theta_b = 0\}$. The first part of the claim follows. The second part follows from the proof of Theorem 3.1. in [13]. $\square$

LEMMA A.2. *Let $\hat{\theta}^{a,\text{ne}_a,\lambda}$ be defined for every $a \in \Gamma(n)$ as in (A.1). Under the assumptions of Theorem 1, for some $c > 0$, for all $a \in \Gamma(n)$,*

$$P(\text{sign}(\hat{\theta}_b^{a,\text{ne}_a,\lambda}) = \text{sign}(\theta_b^a), \forall b \in \text{ne}_a) = 1 - O(\exp(-cn^\varepsilon)) \qquad \text{for } n \to \infty.$$

For the sign-function, it is understood that $\text{sign}(0) = 0$. The lemma says, in other words, that if one could restrict the Lasso estimate to have zero coefficients for all nodes which are not in the neighborhood of node $a$, then the signs of the partial correlations in the neighborhood of node $a$ are estimated consistently under the given assumptions.

PROOF. Using Bonferroni's inequality, and $|ne_a| = o(n)$ for $n \to \infty$, it suffices to show that there exists some $c > 0$ so that for every $a, b \in \Gamma(n)$ with $b \in \text{ne}_a$,

$$P(\text{sign}(\hat{\theta}_b^{a,\text{ne}_a,\lambda}) = \text{sign}(\theta_b^a)) = 1 - O(\exp(-cn^\varepsilon)) \qquad \text{for } n \to \infty.$$

Consider the definition of $\hat{\theta}^{a,\text{ne}_a,\lambda}$ in (A.1),

(A.2) $$\hat{\theta}^{a,\text{ne}_a,\lambda} = \underset{\theta : \theta_k = 0 \, \forall k \notin \text{ne}_a}{\arg\min} (n^{-1}\|\mathbf{X}_a - \mathbf{X}\theta\|_2^2 + \lambda\|\theta\|_1).$$

Assume now that component $b$ of this estimate is fixed at a constant value $\beta$. Denote this new estimate by $\tilde{\theta}^{a,b,\lambda}(\beta)$,

(A.3) $$\tilde{\theta}^{a,b,\lambda}(\beta) = \underset{\theta \in \Theta_{a,b}(\beta)}{\arg\min} (n^{-1}\|\mathbf{X}_a - \mathbf{X}\theta\|_2^2 + \lambda\|\theta\|_1),$$

where

$$\Theta_{a,b}(\beta) := \{\theta \in \mathbb{R}^{p(n)} : \theta_b = \beta; \theta_k = 0, \forall k \notin \text{ne}_a\}.$$

There always exists a value $\beta$ (namely $\beta = \hat{\theta}_b^{a,\text{ne}_a,\lambda}$) so that $\tilde{\theta}^{a,b,\lambda}(\beta)$ is identical to $\hat{\theta}^{a,\text{ne}_a,\lambda}$. Thus, if $\text{sign}(\hat{\theta}_b^{a,\text{ne}_a,\lambda}) \neq \text{sign}(\theta_b^a)$, there would exist some $\beta$ with $\text{sign}(\beta)\text{sign}(\theta_b^a) \leq 0$ so that $\tilde{\theta}^{a,b,\lambda}(\beta)$ would be a solution to (A.2). Using $\text{sign}(\theta_b^a) \neq 0$ for all $b \in \text{ne}_a$, it is thus sufficient to show that for every $\beta$ with $\text{sign}(\beta)\text{sign}(\theta_b^a) < 0$, $\tilde{\theta}^{a,b,\lambda}(\beta)$ cannot be a solution to (A.2) with high probability.

We focus in the following on the case where $\theta_b^a > 0$ for notational simplicity. The case $\theta_b^a < 0$ follows analogously. If $\theta_b^a > 0$, it follows by Lemma A.1 that $\tilde{\theta}^{a,b,\lambda}(\beta)$ with $\tilde{\theta}^{a,b,\lambda}_b(\beta) = \beta \leq 0$ can be a solution to (A.2) only if $G_b(\tilde{\theta}^{a,b,\lambda}(\beta)) \geq -\lambda$. Hence it suffices to show that for some $c > 0$ and all $b \in \text{ne}_a$ with $\theta^a_b > 0$, for $n \to \infty$,

$$(\text{A.4}) \qquad P\left(\sup_{\beta \leq 0}\{G_b(\tilde{\theta}^{a,b,\lambda}(\beta))\} < -\lambda\right) = 1 - O(\exp(-cn^\varepsilon)).$$

Let in the following $\mathbf{R}^\lambda_a(\beta)$ be the $n$-dimensional vector of residuals,

$$(\text{A.5}) \qquad \mathbf{R}^\lambda_a(\beta) := \mathbf{X}_a - \mathbf{X}\tilde{\theta}^{a,b,\lambda}(\beta).$$

We can write $X_b$ as

$$(\text{A.6}) \qquad X_b = \sum_{k \in \text{ne}_a \setminus \{b\}} \theta^{b,\text{ne}_a \setminus \{b\}}_k X_k + W_b,$$

where $W_b$ is independent of $\{X_k; k \in \text{ne}_a \setminus \{b\}\}$. By straightforward calculation, using (A.6),

$$G_b(\tilde{\theta}^{a,b,\lambda}(\beta)) = -2n^{-1}\langle \mathbf{R}^\lambda_a(\beta), \mathbf{W}_b \rangle - \sum_{k \in \text{ne}_a \setminus \{b\}} \theta^{b,\text{ne}_a \setminus \{b\}}_k (2n^{-1}\langle \mathbf{R}^\lambda_a(\beta), \mathbf{X}_k \rangle).$$

By Lemma A.1, for all $k \in \text{ne}_a \setminus \{b\}$, $|G_k(\tilde{\theta}^{a,b,\lambda}(\beta))| = |2n^{-1}\langle \mathbf{R}^\lambda_a(\beta), \mathbf{X}_k \rangle| \leq \lambda$. This together with the equation above yields

$$(\text{A.7}) \qquad G_b(\tilde{\theta}^{a,b,\lambda}(\beta)) \leq -2n^{-1}\langle \mathbf{R}^\lambda_a(\beta), \mathbf{W}_b \rangle + \lambda \|\theta^{b,\text{ne}_a \setminus \{b\}}\|_1.$$

Using Assumption 4, there exists some $\vartheta < \infty$, so that $\|\theta^{b,\text{ne}_a \setminus \{b\}}\|_1 \leq \vartheta$. For proving (A.4) it is therefore sufficient to show that there exists for every $g > 0$ some $c > 0$ so that for all $b \in \text{ne}_a$ with $\theta^a_b > 0$, for $n \to \infty$,

$$(\text{A.8}) \qquad P\left(\inf_{\beta \leq 0}\{2n^{-1}\langle \mathbf{R}^\lambda_a(\beta), \mathbf{W}_b \rangle\} > g\lambda\right) = 1 - O(\exp(-cn^\varepsilon)).$$

With a little abuse of notation, let $W^\| \subseteq \mathbb{R}^n$ be the at most $(|\text{ne}_a| - 1)$-dimensional space which is spanned by the vectors $\{\mathbf{X}_k, k \in \text{ne}_a \setminus \{b\}\}$ and let $W^\perp$ be the orthogonal complement of $W^\|$ in $\mathbb{R}^n$. Split the $n$-dimensional vector $\mathbf{W}_b$ of observations of $W_b$ into the sum of two vectors

$$(\text{A.9}) \qquad \mathbf{W}_b = \mathbf{W}^\perp_b + \mathbf{W}^\|_b,$$

where $\mathbf{W}^\|_b$ is contained in the space $W^\| \subseteq \mathbb{R}^n$, while the remaining part $\mathbf{W}^\perp_b$ is chosen orthogonal to this space (in the orthogonal complement $W^\perp$ of $W^\|$). The inner product in (A.8) can be written as

$$(\text{A.10}) \quad 2n^{-1}\langle \mathbf{R}^\lambda_a(\beta), \mathbf{W}_b \rangle = 2n^{-1}\langle \mathbf{R}^\lambda_a(\beta), \mathbf{W}^\perp_b \rangle + 2n^{-1}\langle \mathbf{R}^\lambda_a(\beta), \mathbf{W}^\|_b \rangle.$$

By Lemma A.3 (see below), there exists for every $g > 0$ some $c > 0$ so that, for $n \to \infty$,

$$P\left(\inf_{\beta \leq 0}\{2n^{-1}\langle \mathbf{R}_a^\lambda(\beta), \mathbf{W}_b^\| \rangle/(1+|\beta|)\} > -g\lambda\right) = 1 - O(\exp(-cn^\varepsilon)).$$

To show (A.8), it is sufficient to prove that there exists for every $g > 0$ some $c > 0$ so that, for $n \to \infty$,

(A.11) $P\left(\inf_{\beta \leq 0}\{2n^{-1}\langle \mathbf{R}_a^\lambda(\beta), \mathbf{W}_b^\perp \rangle - g(1+|\beta|)\lambda\} > g\lambda\right) = 1 - O(\exp(-cn^\varepsilon)).$

For some random variable $V_a$, independent of $X_{\mathrm{ne}_a}$, we have

$$X_a = \sum_{k \in \mathrm{ne}_a} \theta_k^a X_k + V_a.$$

Note that $V_a$ and $W_b$ are independent normally distributed random variables with variances $\sigma_a^2$ and $\sigma_b^2$, respectively. By Assumption 2, $0 < v^2 \leq \sigma_b^2, \sigma_a^2 \leq 1$. Note furthermore that $W_b$ and $X_{\mathrm{ne}_a \setminus \{b\}}$ are independent. Using $\theta^a = \theta^{a,\mathrm{ne}_a}$ and (A.6),

(A.12) $$X_a = \sum_{k \in \mathrm{ne}_a \setminus \{b\}} (\theta_k^a + \theta_b^a \theta_k^{b,\mathrm{ne}_a \setminus \{b\}})X_k + \theta_b^a W_b + V_a.$$

Using (A.12), the definition of the residuals in (A.5) and the orthogonality property of $\mathbf{W}_b^\perp$,

(A.13) $$\begin{aligned}2n^{-1}\langle \mathbf{R}_a^\lambda(\beta), \mathbf{W}_b^\perp \rangle &= 2n^{-1}(\theta_b^a - \beta)\langle \mathbf{W}_b^\perp, \mathbf{W}_b^\perp \rangle + 2n^{-1}\langle \mathbf{V}_a, \mathbf{W}_b^\perp \rangle,\\ &\geq 2n^{-1}(\theta_b^a - \beta)\langle \mathbf{W}_b^\perp, \mathbf{W}_b^\perp \rangle - |2n^{-1}\langle \mathbf{V}_a, \mathbf{W}_b^\perp \rangle|.\end{aligned}$$

The second term, $|2n^{-1}\langle \mathbf{V}_a, \mathbf{W}_b^\perp \rangle|$, is stochastically smaller than $|2n^{-1}\langle \mathbf{V}_a, \mathbf{W}_b \rangle|$ (this can be derived by conditioning on $\{\mathbf{X}_k; k \in \mathrm{ne}_a\}$). Due to independence of $V_a$ and $W_b$, $E(V_a W_b) = 0$. Using Bernstein's inequality (Lemma 2.2.11 in [17]), and $\lambda \sim dn^{-(1-\varepsilon)/2}$ with $\varepsilon > 0$, there exists for every $g > 0$ some $c > 0$ so that

(A.14) $$\begin{aligned}P(|2n^{-1}\langle \mathbf{V}_a, \mathbf{W}_b^\perp \rangle| \geq g\lambda) &\leq P(|2n^{-1}\langle \mathbf{V}_a, \mathbf{W}_b \rangle| \geq g\lambda)\\ &= O(\exp(-cn^\varepsilon)).\end{aligned}$$

Instead of (A.11), it is sufficient by (A.13) and (A.14) to show that there exists for every $g > 0$ a $c > 0$ so that, for $n \to \infty$,

(A.15) $$\begin{aligned}P\Big(\inf_{\beta \leq 0}\{2n^{-1}(\theta_b^a - \beta)\langle \mathbf{W}_b^\perp, \mathbf{W}_b^\perp \rangle &- g(1+|\beta|)\lambda\} > 2g\lambda\Big)\\ &= 1 - O(\exp(-cn^\varepsilon)).\end{aligned}$$

Note that $\sigma_b^{-2}\langle \mathbf{W}_b^\perp, \mathbf{W}_b^\perp \rangle$ follows a $\chi^2_{n-|\mathrm{ne}_a|}$ distribution. As $|\mathrm{ne}_a| = o(n)$ and $\sigma_b^2 \geq v^2$ (by Assumption 2), it follows that there exists some $k > 0$ so that for $n \geq n_0$ with some $n_0(k) \in \mathbb{N}$, and any $c > 0$,

$$P(2n^{-1}\langle \mathbf{W}_b^\perp, \mathbf{W}_b^\perp \rangle > k) = 1 - O(\exp(-cn^\varepsilon)).$$

To show (A.15), it hence suffices to prove that for every $k, \ell > 0$ there exists some $n_0(k, \ell) \in \mathbb{N}$ so that, for all $n \geq n_0$,

(A.16) $$\inf_{\beta \leq 0} \{(\theta_b^a - \beta)k - \ell(1 + |\beta|)\lambda\} > 0.$$

By Assumption 5, $|\pi_{ab}|$ is of order at least $n^{-(1-\xi)/2}$. Using

$$\pi_{ab} = \theta_b^a / (\operatorname{Var}(X_a | X_{\Gamma(n) \setminus \{a\}}) \operatorname{Var}(X_b | X_{\Gamma(n) \setminus \{b\}}))^{1/2}$$

and Assumption 2, this implies that there exists some $q > 0$ so that $\theta_b^a \geq qn^{-(1-\xi)/2}$. As $\lambda \sim dn^{-(1-\varepsilon)/2}$ and, by the assumptions in Theorem 1, $\xi > \varepsilon$, it follows that for every $k, \ell > 0$ and large enough values of $n$,

$$\theta_b^a k - \ell \lambda > 0.$$

It remains to show that for any $k, \ell > 0$ there exists some $n_0(k, \ell)$ so that for all $n \geq n_0$,

$$\inf_{\beta \leq 0} \{-\beta k - \ell |\beta| \lambda\} \geq 0.$$

This follows as $\lambda \to 0$ for $n \to \infty$, which completes the proof. $\square$

LEMMA A.3. *Assume the conditions of Theorem 1 hold. Let $\mathbf{R}_a^\lambda(\beta)$ be defined as in (A.5) and $\mathbf{W}_b^\|$ as in (A.9). For any $g > 0$ there exists $c > 0$ so that for all $a, b \in \Gamma(n)$, for $n \to \infty$,*

$$P\left(\sup_{\beta \in \mathbb{R}} |2n^{-1} \langle \mathbf{R}_a^\lambda(\beta), \mathbf{W}_b^\| \rangle| / (1 + |\beta|) < g\lambda\right) = 1 - O(\exp(-cn^\varepsilon)).$$

PROOF. By Schwarz's inequality,

(A.17) $$|2n^{-1} \langle \mathbf{R}_a^\lambda(\beta), \mathbf{W}_b^\| \rangle| / (1 + |\beta|) \leq 2n^{-1/2} \|\mathbf{W}_b^\|\|_2 \frac{n^{-1/2} \|\mathbf{R}_a^\lambda(\beta)\|_2}{1 + |\beta|}.$$

The sum of squares of the residuals is increasing with increasing value of $\lambda$. Thus, $\|\mathbf{R}_a^\lambda(\beta)\|_2^2 \leq \|\mathbf{R}_a^\infty(\beta)\|_2^2$. By definition of $\mathbf{R}_a^\lambda$ in (A.5), and using (A.3),

$$\|\mathbf{R}_a^\infty(\beta)\|_2^2 = \|\mathbf{X}_a - \beta \mathbf{X}_b\|_2^2,$$

and hence

$$\|\mathbf{R}_a^\lambda(\beta)\|_2^2 \leq (1 + |\beta|)^2 \max\{\|\mathbf{X}_a\|_2^2, \|\mathbf{X}_b\|_2^2\}.$$

Therefore, for any $q > 0$,

$$P\left(\sup_{\beta \in \mathbb{R}} \frac{n^{-1/2} \|\mathbf{R}_a^\lambda(\beta)\|_2}{1 + |\beta|} > q\right) \leq P(n^{-1/2} \max\{\|\mathbf{X}_a\|_2, \|\mathbf{X}_b\|_2\} > q).$$

Note that both $\|\mathbf{X}_a\|_2^2$ and $\|\mathbf{X}_b\|_2^2$ are $\chi_n^2$-distributed. Thus there exist $q > 1$ and $c > 0$ so that

$$\text{(A.18)} \quad P\left(\sup_{\beta \in \mathbb{R}} \frac{n^{-1/2}\|\mathbf{R}_a^\lambda(\beta)\|_2}{1 + |\beta|} > q\right) = O(\exp(-cn^\varepsilon)) \quad \text{for } n \to \infty.$$

It remains to show that for every $g > 0$ there exists some $c > 0$ so that

$$\text{(A.19)} \qquad P(n^{-1/2}\|\mathbf{W}_b^\|\|_2 > g\lambda) = O(\exp(-cn^\varepsilon)) \quad \text{for } n \to \infty.$$

The expression $\sigma_b^{-2}\langle \mathbf{W}_b^\|, \mathbf{W}_b^\| \rangle$ is $\chi_{|\text{ne}_a|-1}^2$-distributed. As $\sigma_b \le 1$ and $|\text{ne}_a| = O(n^\kappa)$, it follows that $n^{-1/2}\|\mathbf{W}_b^\|\|_2$ is for some $t > 0$ stochastically smaller than

$$tn^{-(1-\kappa)/2}(Z/n^\kappa)^{1/2},$$

where $Z$ is $\chi_{n^\kappa}^2$-distributed. Thus, for every $g > 0$,

$$P(n^{-1/2}\|\mathbf{W}_b^\|\|_2 > g\lambda) \le P((Z/n^\kappa) > (g/t)^2 n^{(1-\kappa)}\lambda^2).$$

As $\lambda^{-1} = O(n^{(1-\varepsilon)/2})$, it follows that $n^{1-\kappa}\lambda^2 \ge hn^{\varepsilon-\kappa}$ for some $h > 0$ and sufficiently large $n$. By the properties of the $\chi^2$ distribution and $\varepsilon > \kappa$, by assumption in Theorem 1, claim (A.19) follows. This completes the proof. □

PROOF OF PROPOSITION 1. All diagonal elements of the covariance matrices $\Sigma(n)$ are equal to 1, while all off-diagonal elements vanish for all pairs except for $a, b \in \Gamma(n)$, where $\Sigma_{ab}(n) = s$ with $0 < s < 1$. Assume w.l.o.g. that $a$ corresponds to the first and $b$ to the second variable. The best vector of coefficients $\theta^a$ for linear prediction of $X_a$ is given by $\theta^a = (0, -K_{ab}/K_{aa}, 0, 0, \ldots) = (0, s, 0, 0, \ldots)$, where $K = \Sigma^{-1}(n)$. A necessary condition for $\tilde{\text{ne}}_a^\lambda = \text{ne}_a$ is that $\hat{\theta}^{a,\lambda} = (0, \tau, 0, 0, \ldots)$ is the oracle Lasso solution for some $\tau \ne 0$. In the following, we show first that

$$\text{(A.20)} \qquad P(\exists \lambda, \tau \ge s : \hat{\theta}^{a,\lambda} = (0, \tau, 0, 0, \ldots)) \to 0, \qquad n \to \infty.$$

The proof is then completed by showing in addition that $(0, \tau, 0, 0, \ldots)$ cannot be the *oracle* Lasso solution as long as $\tau < s$.

We begin by showing (A.20). If $\hat{\theta} = (0, \tau, 0, 0, \ldots)$ is a Lasso solution for some value of the penalty, it follows that, using Lemma A.1 and positivity of $\tau$,

$$\text{(A.21)} \quad \langle \mathbf{X}_1 - \tau\mathbf{X}_2, \mathbf{X}_2 \rangle \ge |\langle \mathbf{X}_1 - \tau\mathbf{X}_2, \mathbf{X}_k \rangle| \qquad \forall k \in \Gamma(n), \ k > 2.$$

Under the given assumptions, $X_2, X_3, \ldots$ can be understood to be independently and identically distributed, while $X_1 = sX_2 + W_1$, with $W_1$ independent of $(X_2, X_3, \ldots)$. Substituting $X_1 = sX_2 + W_1$ in (A.21) yields for all $k \in \Gamma(n)$ with $k > 2$,

$$\langle \mathbf{W}_1, \mathbf{X}_2 \rangle - (\tau - s)\langle \mathbf{X}_2, \mathbf{X}_2 \rangle \ge |\langle \mathbf{W}_1, \mathbf{X}_k \rangle - (\tau - s)\langle \mathbf{X}_2, \mathbf{X}_k \rangle|.$$

Let $U_2, U_3, \ldots, U_{p(n)}$ be the random variables defined by $U_k = \langle \mathbf{W}_1, \mathbf{X}_k \rangle$. Note that the random variables $U_k$, $k = 2, \ldots, p(n)$, are exchangeable. Let furthermore

$$D = \langle \mathbf{X}_2, \mathbf{X}_2 \rangle - \max_{k \in \Gamma(n), k > 2} |\langle \mathbf{X}_2, \mathbf{X}_k \rangle|.$$

The inequality above implies then

$$U_2 > \max_{k \in \Gamma(n), k > 2} U_k + (\tau - s) D.$$

To show the claim, it thus suffices to show that

$$(A.22) \qquad P\left(U_2 > \max_{k \in \Gamma(n), k > 2} U_k + (\tau - s) D\right) \to 0 \qquad \text{for } n \to \infty.$$

Using $\tau - s > 0$,

$$P\left(U_2 > \max_{k \in \Gamma(n), k > 2} U_k + (\tau - s) D\right) \leq P\left(U_2 > \max_{k \in \Gamma(n), k > 2} U_k\right) + P(D < 0).$$

Using the assumption that $s < 1$, it follows by $p(n) = o(n^\gamma)$ for some $\gamma > 0$ and a Bernstein-type inequality that

$$P(D < 0) \to 0 \qquad \text{for } n \to \infty.$$

Furthermore, as $U_2, \ldots, U_{p(n)}$ are exchangeable,

$$P\left(U_2 > \max_{k \in \Gamma(n), k > 2} U_k\right) = (p(n) - 1)^{-1} \to 0 \qquad \text{for } n \to \infty,$$

which shows that (A.22) holds. The claim (A.20) follows.

It hence suffices to show that $(0, \tau, 0, 0, \ldots)$ with $\tau < s$ cannot be the *oracle* Lasso solution. Let $\tau_{\max}$ be the maximal value of $\tau$ so that $(0, \tau, 0, \ldots)$ is a Lasso solution for some value $\lambda > 0$. By the previous assumption, $\tau_{\max} < s$. For $\tau < \tau_{\max}$, the vector $(0, \tau, 0, \ldots)$ cannot be the oracle Lasso solution. We show in the following that $(0, \tau_{\max}, 0, \ldots)$ cannot be an oracle Lasso solution either. Suppose that $(0, \tau_{\max}, 0, 0, \ldots)$ is the Lasso solution $\hat{\theta}^{a,\lambda}$ for some $\lambda = \tilde{\lambda} > 0$. As $\tau_{\max}$ is the maximal value such that $(0, \tau, 0, \ldots)$ is a Lasso solution, there exists some $k \in \Gamma(n) > 2$, such that

$$|n^{-1} \langle \mathbf{X}_1 - \tau_{\max} \mathbf{X}_2, \mathbf{X}_2 \rangle| = |n^{-1} \langle \mathbf{X}_1 - \tau_{\max} \mathbf{X}_2, \mathbf{X}_k \rangle|,$$

and the value of both components $G_2$ and $G_k$ of the gradient is equal to $\tilde{\lambda}$. By appropriately reordering the variables we can assume that $k = 3$. Furthermore, it holds a.s. that

$$\max_{k \in \Gamma(n), k > 3} |\langle \mathbf{X}_1 - \tau_{\max} \mathbf{X}_2, \mathbf{X}_k \rangle| < \tilde{\lambda}.$$

Hence, for sufficiently small $\delta\lambda \geq 0$, a Lasso solution for the penalty $\tilde{\lambda} - \delta\lambda$ is given by

$$(0, \tau_{\max} + \delta\theta_2, \delta\theta_3, 0, \ldots).$$

Let $H_n$ be the empirical covariance matrix of $(X_2, X_3)$. Assume w.l.o.g. that $n^{-1}\langle \mathbf{X}_1 - \tau_{\max}\mathbf{X}_2, \mathbf{X}_k \rangle > 0$ and $n^{-1}\langle \mathbf{X}_2, \mathbf{X}_2 \rangle = n^{-1}\langle \mathbf{X}_3, \mathbf{X}_3 \rangle = 1$. Following, for example, Efron et al. ([6], page 417), the components $(\delta\theta_2, \delta\theta_3)$ are then given by $H_n^{-1}(1,1)^T$, from which it follows that $\delta\theta_2 = \delta\theta_3$, which we abbreviate by $\delta\theta$ in the following (one can accommodate a negative sign for $n^{-1}\langle \mathbf{X}_1 - \tau_{\max}\mathbf{X}_2, \mathbf{X}_k \rangle$ by reversing the sign of $\delta\theta_3$). Denote by $L_\delta$ the squared error loss for this solution. Then, for sufficiently small $\delta\theta$,

$$L_\delta - L_0 = E(X_1 - (\tau_{\max} + \delta\theta)X_2 + \delta\theta X_3)^2 - E(X_1 - \tau_{\max}X_2)^2$$
$$= (s - (\tau_{\max} + \delta\theta))^2 + \delta\theta^2 - (s - \tau_{\max})^2$$
$$= -2(s - \tau_{\max})\delta\theta + 2\delta\theta^2.$$

It holds that $L_{\delta\theta} - L_0 < 0$ for any $0 < \delta\theta < 1/2(s - \tau_{\max})$, which shows that $(0, \tau, 0, \ldots)$ cannot be the oracle solution for $\tau < s$. Together with (A.20), this completes the proof. □

PROOF OF PROPOSITION 2. The subdifferential of the argument in (6),

$$E\left(X_a - \sum_{m \in \Gamma(n)} \theta_m^a X_m\right)^2 + \eta\|\theta^a\|_1,$$

with respect to $\theta_k^a$, $k \in \Gamma(n) \setminus \{a\}$, is given by

$$-2E\left(\left(X_a - \sum_{m \in \Gamma(n)} \theta_m^a X_m\right)X_k\right) + \eta e_k,$$

where $e_k \in [-1, 1]$ if $\theta_k^a = 0$, and $e_k = \text{sign}(\theta_k^a)$ if $\theta_k^a \neq 0$. Using the fact that $\text{ne}_a(\eta) = \text{ne}_a$, it follows as in Lemma A.1 that for all $k \in \text{ne}_a$,

(A.23) $$2E\left(\left(X_a - \sum_{m \in \Gamma(n)} \theta_m^a(\eta) X_m\right)X_k\right) = \eta \,\text{sign}(\theta_k^a)$$

and, for $b \notin \text{ne}_a$,

(A.24) $$\left|2E\left(\left(X_a - \sum_{m \in \Gamma(n)} \theta_m^a(\eta) X_m\right)X_b\right)\right| \leq \eta.$$

A variable $X_b$ with $b \notin \text{ne}_a$ can be written as

$$X_b = \sum_{k \in \text{ne}_a} \theta_k^{b, \text{ne}_a} X_k + W_b,$$

where $W_b$ is independent of $\{X_k; k \in \mathrm{cl}_a\}$. Using this in (A.24) yields

$$\left| 2 \sum_{k \in \mathrm{ne}_a} \theta_k^{b,\mathrm{ne}_a} E\left( \left( X_a - \sum_{m \in \Gamma(n)} \theta_m^a(\eta) X_m \right) X_k \right) \right| \leq \eta.$$

Using (A.23) and $\theta^a = \theta^{a,\mathrm{ne}_a}$, it follows that

$$\left| \sum_{k \in \mathrm{ne}_a} \theta_k^{b,\mathrm{ne}_a} \mathrm{sign}(\theta_k^{a,\mathrm{ne}_a}) \right| \leq 1,$$

which completes the proof. $\square$

PROOF OF THEOREM 1. The event $\hat{\mathrm{ne}}_a^\lambda \not\subseteq \mathrm{ne}_a$ is equivalent to the event that there exists some node $b \in \Gamma(n) \setminus \mathrm{cl}_a$ in the set of nonneighbors of node $a$ such that the estimated coefficient $\hat{\theta}_b^{a,\lambda}$ is not zero. Thus

(A.25) $\qquad P(\hat{\mathrm{ne}}_a^\lambda \subseteq \mathrm{ne}_a) = 1 - P(\exists b \in \Gamma(n) \setminus \mathrm{cl}_a : \hat{\theta}_b^{a,\lambda} \neq 0).$

Consider the Lasso estimate $\hat{\theta}^{a,\mathrm{ne}_a,\lambda}$, which is by (A.1) constrained to have nonzero components only in the neighborhood of node $a \in \Gamma(n)$. Using $|\mathrm{ne}_a| = O(n^\kappa)$ with some $\kappa < 1$, we can assume w.l.o.g. that $|\mathrm{ne}_a| \leq n$. This in turn implies (see, e.g., [13]) that $\hat{\theta}^{a,\mathrm{ne}_a,\lambda}$ is a.s. a unique solution to (A.1) with $\mathcal{A} = \mathrm{ne}_a$. Let $\mathcal{E}$ be the event

$$\max_{k \in \Gamma(n) \setminus \mathrm{cl}_a} |G_k(\hat{\theta}^{a,\mathrm{ne}_a,\lambda})| < \lambda.$$

Conditional on the event $\mathcal{E}$, it follows from the first part of Lemma A.1 that $\hat{\theta}^{a,\mathrm{ne}_a,\lambda}$ is not only a solution of (A.1), with $\mathcal{A} = \mathrm{ne}_a$, but as well a solution of (3), where $\mathcal{A} = \Gamma(n) \setminus \{a\}$. As $\hat{\theta}_b^{a,\mathrm{ne}_a,\lambda} = 0$ for all $b \in \Gamma(n) \setminus \mathrm{cl}_a$, it follows from the second part of Lemma A.1 that $\hat{\theta}_b^{a,\lambda} = 0, \forall b \in \Gamma(n) \setminus \mathrm{cl}_a$. Hence

$$P(\exists b \in \Gamma(n) \setminus \mathrm{cl}_a : \hat{\theta}_b^{a,\lambda} \neq 0) \leq 1 - P(\mathcal{E})$$

$$= P\left( \max_{k \in \Gamma(n) \setminus \mathrm{cl}_a} |G_k(\hat{\theta}^{a,\mathrm{ne}_a,\lambda})| \geq \lambda \right),$$

where

(A.26) $\qquad G_b(\hat{\theta}^{a,\mathrm{ne}_a,\lambda}) = -2n^{-1}\langle \mathbf{X}_a - \mathbf{X}\hat{\theta}^{a,\mathrm{ne}_a,\lambda}, \mathbf{X}_b \rangle.$

Using Bonferroni's inequality and $p(n) = O(n^\gamma)$ for any $\gamma > 0$, it suffices to show that there exists a constant $c > 0$ so that for all $b \in \Gamma(n) \setminus \mathrm{cl}_a$,

(A.27) $\qquad P(|G_b(\hat{\theta}^{a,\mathrm{ne}_a,\lambda})| \geq \lambda) = O(\exp(-cn^\varepsilon)).$

One can write for any $b \in \Gamma(n) \setminus \mathrm{cl}_a$,

(A.28) $\qquad X_b = \sum_{m \in \mathrm{ne}_a} \theta_m^{b,\mathrm{ne}_a} X_m + V_b,$

where $V_b \sim \mathcal{N}(0, \sigma_b^2)$ for some $\sigma_b^2 \leq 1$ and $V_b$ is independent of $\{X_m; m \in \mathrm{cl}_a\}$. Hence

$$G_b(\hat{\theta}^{a,\mathrm{ne}_a,\lambda}) = -2n^{-1} \sum_{m \in \mathrm{ne}_a} \theta_m^{b,\mathrm{ne}_a} \langle \mathbf{X}_a - \mathbf{X}\hat{\theta}^{a,\mathrm{ne}_a,\lambda}, \mathbf{X}_m \rangle$$

$$- 2n^{-1} \langle \mathbf{X}_a - \mathbf{X}\hat{\theta}^{a,\mathrm{ne}_a,\lambda}, \mathbf{V}_b \rangle.$$

By Lemma A.2, there exists some $c > 0$ so that with probability $1 - O(\exp(-cn^\varepsilon))$,

(A.29) $$\mathrm{sign}(\hat{\theta}_k^{a,\mathrm{ne}_a,\lambda}) = \mathrm{sign}(\theta_k^{a,\mathrm{ne}_a}) \qquad \forall k \in \mathrm{ne}_a.$$

In this case by Lemma A.1

$$2n^{-1} \sum_{m \in \mathrm{ne}_a} \theta_m^{b,\mathrm{ne}_a} \langle \mathbf{X}_a - \mathbf{X}\hat{\theta}^{a,\mathrm{ne}_a,\lambda}, \mathbf{X}_m \rangle = \left( \sum_{m \in \mathrm{ne}_a} \mathrm{sign}(\theta_m^{a,\mathrm{ne}_a}) \theta_m^{b,\mathrm{ne}_a} \right) \lambda.$$

If (A.29) holds, the gradient is given by

(A.30) $$G_b(\hat{\theta}^{a,\mathrm{ne}_a,\lambda}) = -\left( \sum_{m \in \mathrm{ne}_a} \mathrm{sign}(\theta_m^{a,\mathrm{ne}_a}) \theta_m^{b,\mathrm{ne}_a} \right) \lambda - 2n^{-1} \langle \mathbf{X}_a - \mathbf{X}\hat{\theta}^{a,\mathrm{ne}_a,\lambda}, \mathbf{V}_b \rangle.$$

Using Assumption 6 and Proposition 2, there exists some $\delta < 1$ so that

$$\left| \sum_{m \in \mathrm{ne}_a} \mathrm{sign}(\theta_m^{a,\mathrm{ne}_a}) \theta_m^{b,\mathrm{ne}_a} \right| \leq \delta.$$

The absolute value of the coefficient $G_b$ of the gradient in (A.26) is hence bounded with probability $1 - O(\exp(-cn^\varepsilon))$ by

(A.31) $$|G_b(\hat{\theta}^{a,\mathrm{ne}_a,\lambda})| \leq \delta\lambda + |2n^{-1} \langle \mathbf{X}_a - \mathbf{X}\hat{\theta}^{a,\mathrm{ne}_a,\lambda}, \mathbf{V}_b \rangle|.$$

Conditional on $\mathbf{X}_{\mathrm{cl}_a} = \{\mathbf{X}_k; k \in \mathrm{cl}_a\}$, the random variable

$$\langle \mathbf{X}_a - \mathbf{X}\hat{\theta}_k^{a,\mathrm{ne}_a,\lambda}, \mathbf{V}_b \rangle$$

is normally distributed with mean zero and variance $\sigma_b^2 \|\mathbf{X}_a - \mathbf{X}\hat{\theta}^{a,\mathrm{ne}_a,\lambda}\|_2^2$. On the one hand, $\sigma_b^2 \leq 1$. On the other hand, by definition of $\hat{\theta}^{a,\mathrm{ne}_a,\lambda}$,

$$\|\mathbf{X}_a - \mathbf{X}\hat{\theta}^{a,\mathrm{ne}_a,\lambda}\|_2 \leq \|\mathbf{X}_a\|_2.$$

Thus

$$|2n^{-1} \langle \mathbf{X}_a - \mathbf{X}\hat{\theta}^{a,\mathrm{ne}_a,\lambda}, \mathbf{V}_b \rangle|$$

is stochastically smaller than or equal to $|2n^{-1} \langle \mathbf{X}_a, \mathbf{V}_b \rangle|$. Using (A.31), it remains to be shown that for some $c > 0$ and $\delta < 1$,

$$P(|2n^{-1} \langle \mathbf{X}_a, \mathbf{V}_b \rangle| \geq (1-\delta)\lambda) = O(\exp(-cn^\varepsilon)).$$

As $V_b$ and $X_a$ are independent, $E(X_a V_b) = 0$. Using the Gaussianity and bounded variance of both $X_a$ and $V_b$, there exists some $g < \infty$ so that $E(\exp(|X_a V_b|)) \leq g$. Hence, using Bernstein's inequality and the boundedness of $\lambda$, for some $c > 0$, for all $b \in \text{ne}_a$, $P(|2n^{-1}\langle \mathbf{X}_a, \mathbf{V}_b \rangle| \geq (1-\delta)\lambda) = O(\exp(-cn\lambda^2))$. The claim (A.27) follows, which completes the proof. $\square$

PROOF OF PROPOSITION 3. Following a similar argument as in Theorem 1 up to (A.27), it is sufficient to show that for every $a, b$ with $b \in \Gamma(n) \setminus \text{cl}_a$ and $|S_a(b)| > 1$,

(A.32) $$P(|G_b(\hat{\theta}^{a,\text{ne}_a,\lambda})| > \lambda) \to 1 \qquad \text{for } n \to \infty.$$

Using (A.30) in the proof of Theorem 1, one can conclude that for some $\delta > 1$, with probability converging to 1 for $n \to \infty$,

(A.33) $$|G_b(\hat{\theta}^{a,\text{ne}_a,\lambda})| \geq \delta\lambda - |2n^{-1}\langle \mathbf{X}_a - \mathbf{X}\hat{\theta}^{a,\text{ne}_a,\lambda}, \mathbf{V}_b \rangle|.$$

Using the identical argument as in the proof of Theorem 1 below (A.31), for the second term, for any $g > 0$,

$$P(|2n^{-1}\langle \mathbf{X}_a - \mathbf{X}\hat{\theta}^{a,\text{ne}_a,\lambda}, \mathbf{V}_b \rangle| > g\lambda) \to 0 \qquad \text{for } n \to \infty,$$

which together with $\delta > 1$ in (A.33) shows that (A.32) holds. This completes the proof. $\square$

PROOF OF THEOREM 2. First, $P(\text{ne}_a \subseteq \hat{\text{ne}}_a^\lambda) = 1 - P(\exists b \in \text{ne}_a : \hat{\theta}_b^{a,\lambda} = 0)$. Let $\mathcal{E}$ be again the event

(A.34) $$\max_{k \in \Gamma(n) \setminus \text{cl}_a} |G_k(\hat{\theta}^{a,\text{ne}_a,\lambda})| < \lambda.$$

Conditional on $\mathcal{E}$, we can conclude as in the proof of Theorem 1 that $\hat{\theta}^{a,\text{ne}_a,\lambda}$ and $\hat{\theta}^{a,\lambda}$ are unique solutions to (A.1) and (3), respectively, and $\hat{\theta}^{a,\text{ne}_a,\lambda} = \hat{\theta}^{a,\lambda}$. Thus

$$P(\exists b \in \text{ne}_a : \hat{\theta}_b^{a,\lambda} = 0) \leq P(\exists b \in \text{ne}_a : \hat{\theta}_b^{a,\text{ne}_a,\lambda} = 0) + P(\mathcal{E}^c).$$

It follows from the proof of Theorem 1 that there exists some $c > 0$ so that $P(\mathcal{E}^c) = O(\exp(-cn^\varepsilon))$. Using Bonferroni's inequality, it hence remains to show that there exists some $c > 0$ so that for all $b \in \text{ne}_a$,

(A.35) $$P(\hat{\theta}_b^{a,\text{ne}_a,\lambda} = 0) = O(\exp(-cn^\varepsilon)).$$

This follows from Lemma A.2, which completes the proof. $\square$

PROOF OF PROPOSITION 4. The proof of Proposition 4 is to a large extent analogous to the proofs of Theorems 1 and 2. Let $\mathcal{E}$ be again the event (A.34). Conditional on the event $\mathcal{E}$, we can conclude as before that $\hat{\theta}^{a,\mathrm{ne}_a,\lambda}$ and $\hat{\theta}^{a,\lambda}$ are unique solutions to (A.1) and (3), respectively, and $\hat{\theta}^{a,\mathrm{ne}_a,\lambda} = \hat{\theta}^{a,\lambda}$. Thus, for any $b \in \mathrm{ne}_a$,

$$P(b \notin \hat{\mathrm{ne}}_a^\lambda) = P(\hat{\theta}_b^{a,\lambda} = 0) \geq P(\hat{\theta}_b^{a,\mathrm{ne}_a,\lambda} = 0 | \mathcal{E}) P(\mathcal{E}).$$

Since $P(\mathcal{E}) \to 1$ for $n \to \infty$ by Theorem 1,

$$P(\hat{\theta}_b^{a,\mathrm{ne}_a,\lambda} = 0 | \mathcal{E}) P(\mathcal{E}) \to P(\hat{\theta}_b^{a,\mathrm{ne}_a,\lambda} = 0) \qquad \text{for } n \to \infty.$$

It thus suffices to show that for all $b \in \mathrm{ne}_a$ with $|\pi_{ab}| = O(n^{-(1-\xi)/2})$ and $\xi < \varepsilon$,

$$P(\hat{\theta}_b^{a,\mathrm{ne}_a,\lambda} = 0) \to 1 \qquad \text{for } n \to \infty.$$

This holds if

(A.36) $\qquad P(|G_b(\hat{\theta}^{a,\mathrm{ne}_a \setminus \{b\},\lambda})| < \lambda) \to 1 \qquad \text{for } n \to \infty,$

as $|G_b(\hat{\theta}^{a,\mathrm{ne}_a \setminus \{b\},\lambda})| < \lambda$ implies that $\hat{\theta}^{a,\mathrm{ne}_a \setminus \{b\},\lambda} = \hat{\theta}^{a,\mathrm{ne}_a,\lambda}$ and hence $\hat{\theta}_b^{a,\mathrm{ne}_a,\lambda} = 0$. Using (A.7),

$$|G_b(\hat{\theta}^{a,\mathrm{ne}_a \setminus \{b\},\lambda})| \leq |2n^{-1} \langle \mathbf{R}_a^\lambda(0), \mathbf{W}_b \rangle| + \lambda \|\theta^{b,\mathrm{ne}_a \setminus \{b\}}\|_1.$$

By assumption $\|\theta^{b,\mathrm{ne}_a \setminus \{b\}}\|_1 < 1$. It is thus sufficient to show that for any $g > 0$,

(A.37) $\qquad P(|2n^{-1} \langle \mathbf{R}_a^\lambda(0), \mathbf{W}_b \rangle| < g\lambda) \to 1 \qquad \text{for } n \to \infty.$

Analogously to (A.10), we can write

(A.38) $|2n^{-1} \langle \mathbf{R}_a^\lambda(0), \mathbf{W}_b \rangle| \leq |2n^{-1} \langle \mathbf{R}_a^\lambda(0), \mathbf{W}_b^\perp \rangle| + |2n^{-1} \langle \mathbf{R}_a^\lambda(0), \mathbf{W}_b^\| \rangle|.$

Using Lemma A.3, it follows for the last term on the right-hand side that for every $g > 0$,

$$P(|2n^{-1} \langle \mathbf{R}_a^\lambda(0), \mathbf{W}_b^\| \rangle| < g\lambda) \to 1 \qquad \text{for } n \to \infty.$$

Using (A.13) and (A.14), it hence remains to show that for every $g > 0$,

(A.39) $\qquad P(|2n^{-1} \theta_b^{a,\mathrm{ne}_a} \langle \mathbf{W}_b^\perp, \mathbf{W}_b^\perp \rangle| < g\lambda) \to 1 \qquad \text{for } n \to \infty.$

We have already noted above that the term $\sigma_b^{-2} \langle \mathbf{W}_b^\perp, \mathbf{W}_b^\perp \rangle$ follows a $\chi^2_{n-|\mathrm{ne}_a|}$ distribution and is hence stochastically smaller than a $\chi^2_n$-distributed random variable. By Assumption 2, $\sigma_b^2 \leq 1$. Furthermore, using Assumption 2 and $|\pi_{ab}| = O(n^{-(1-\xi)/2})$, $|\theta_b^{a,\mathrm{ne}_a}| = O(n^{-(1-\xi)/2})$. Hence, with $\lambda \sim dn^{-(1-\varepsilon)/2}$, it follows that for some constant $k > 0$, $\lambda/|\theta_b^{a,\mathrm{ne}_a}| \geq kn^{(\varepsilon-\xi)/2}$. Thus, for some constant $c > 0$,

(A.40) $\qquad P(|2n^{-1} \theta_b^{a,\mathrm{ne}_a} \langle \mathbf{W}_b^\perp, \mathbf{W}_b^\perp \rangle| < g\lambda) \geq P(Z/n < cn^{(\varepsilon-\xi)/2}),$

where $Z$ follows a $\chi_n^2$ distribution. By the properties of the $\chi^2$ distribution and the assumption $\xi < \varepsilon$, the right-hand side in (A.40) converges to 1 for $n \to \infty$, from which (A.39) and hence the claim follow. $\square$

PROOF OF THEOREM 3. A necessary condition for $\hat{C}_a^\lambda \nsubseteq C_a$ is that there exists an edge in $\hat{E}^\lambda$ joining two nodes in two different connectivity components. Hence

$$P(\exists a \in \Gamma(n) : \hat{C}_a^\lambda \nsubseteq C_a) \leq p(n) \max_{a \in \Gamma(n)} P(\exists b \in \Gamma(n) \setminus C_a : b \in \hat{\text{ne}}_a^\lambda).$$

Using the same arguments as in the proof of Theorem 1,

$$P(\exists b \in \Gamma(n) \setminus C_a : b \in \hat{\text{ne}}_a^\lambda) \leq P\left(\max_{b \in \Gamma(n) \setminus C_a} |G_b(\hat{\theta}^{a,C_a,\lambda})| \geq \lambda\right),$$

where $\hat{\theta}^{a,C_a,\lambda}$, according to (A.1), has nonzero components only for variables in the connectivity component $C_a$ of node $a$. Hence it is sufficient to show that

(A.41) $\qquad p(n)^2 \max_{a \in \Gamma(n), b \in \Gamma(n) \setminus C_a} P(|G_b(\hat{\theta}^{a,C_a,\lambda})| \geq \lambda) \leq \alpha.$

The gradient is given by $G_b(\hat{\theta}^{a,C_a,\lambda}) = -2n^{-1}\langle \mathbf{X}_a - \mathbf{X}\hat{\theta}^{a,C_a,\lambda}, \mathbf{X}_b \rangle$. For all $k \in C_a$ the variables $X_b$ and $X_k$ are independent as they are in different connectivity components. Hence, conditional on $\mathbf{X}_{C_a} = \{\mathbf{X}_k; k \in C_a\}$,

$$G_b(\hat{\theta}^{a,C_a,\lambda}) \sim \mathcal{N}(0, R^2/n),$$

where $R^2 = 4n^{-1}\|\mathbf{X}_a - \mathbf{X}\hat{\theta}^{a,C_a,\lambda}\|_2^2$, which is smaller than or equal to $\hat{\sigma}_a^2 = 4n^{-1}\|\mathbf{X}_a\|_2^2$ by definition of $\hat{\theta}^{a,C_a,\lambda}$. Hence for all $a \in \Gamma(n)$ and $b \in \Gamma(n) \setminus C_a$,

$$P(|G_b(\hat{\theta}^{a,C_a,\lambda})| \geq \lambda | \mathbf{X}_{C_a}) \leq 2\tilde{\Phi}(\sqrt{n}\lambda/(2\hat{\sigma}_a)),$$

where $\tilde{\Phi} = 1 - \Phi$. It follows for the $\lambda$ proposed in (9) that $P(|G_b(\hat{\theta}^{a,C_a,\lambda})| \geq \lambda | \mathbf{X}_{C_a}) \leq \alpha p(n)^{-2}$, and therefore $P(|G_b(\hat{\theta}^{a,C_a,\lambda})| \geq \lambda) \leq \alpha p(n)^{-2}$. Thus (A.41) follows, which completes the proof. $\square$

## REFERENCES

[1] BUHL, S. (1993). On the existence of maximum-likelihood estimators for graphical Gaussian models. *Scand. J. Statist.* **20** 263–270. MR1241392
[2] CHEN, S., DONOHO, D. and SAUNDERS, M. (2001). Atomic decomposition by basis pursuit. *SIAM Rev.* **43** 129–159. MR1854649
[3] DEMPSTER, A. (1972). Covariance selection. *Biometrics* **28** 157–175.
[4] DRTON, M. and PERLMAN, M. (2004). Model selection for Gaussian concentration graphs. *Biometrika* **91** 591–602. MR2090624
[5] EDWARDS, D. (2000). *Introduction to Graphical Modelling*, 2nd ed. Springer, New York. MR1880319

[6] EFRON, B., HASTIE, T., JOHNSTONE, I. and TIBSHIRANI, R. (2004). Least angle regression (with discussion). *Ann. Statist.* **32** 407–499. MR2060166
[7] FRANK, I. and FRIEDMAN, J. (1993). A statistical view of some chemometrics regression tools (with discussion). *Technometrics* **35** 109–148.
[8] GREENSHTEIN, E. and RITOV, Y. (2004). Persistence in high-dimensional linear predictor selection and the virtue of over-parametrization. *Bernoulli* **10** 971–988. MR2108039
[9] HECKERMAN, D., CHICKERING, D. M., MEEK, C., ROUNTHWAITE, R. and KADIE, C. (2000). Dependency networks for inference, collaborative filtering and data visualization. *J. Machine Learning Research* **1** 49–75.
[10] JUDITSKY, A. and NEMIROVSKI, A. (2000). Functional aggregation for nonparametric regression. *Ann. Statist.* **28** 681–712. MR1792783
[11] KNIGHT, K. and FU, W. (2000). Asymptotics for lasso-type estimators. *Ann. Statist.* **28** 1356–1378. MR1805787
[12] LAURITZEN, S. (1996). *Graphical Models.* Clarendon Press, Oxford. MR1419991
[13] OSBORNE, M., PRESNELL, B. and TURLACH, B. (2000). On the lasso and its dual. *J. Comput. Graph. Statist.* **9** 319–337. MR1822089
[14] SHAO, J. (1993). Linear model selection by cross-validation. *J. Amer. Statist. Assoc.* **88** 486–494. MR1224373
[15] SPEED, T. and KIIVERI, H. (1986). Gaussian Markov distributions over finite graphs. *Ann. Statist.* **14** 138–150. MR0829559
[16] TIBSHIRANI, R. (1996). Regression shrinkage and selection via the lasso. *J. Roy. Statist. Soc. Ser. B* **58** 267–288. MR1379242
[17] VAN DER VAART, A. and WELLNER, J. (1996). *Weak Convergence and Empirical Processes*: *With Applications to Statistics.* Springer, New York. MR1385671

SEMINAR FÜR STATISTIK
ETH ZÜRICH
CH-8092 ZÜRICH
SWITZERLAND
E-MAIL: nicolai@stat.math.ethz.ch
buhlmann@stat.math.ethz.ch